%% file: agt-3-24.tex
\documentclass{gtart}

\input agtout

\lognumber{24}
\volumenumber{3}
\volumeyear{2003}
\papernumber{24}
\published{30 July 2003}
\pagenumbers{709}{718}
\received{7 January 2003}
\accepted{14 July 2003}

\usepackage{amssymb,amsmath}%,amscd} 
   
\newtheorem{theorem}{Theorem}[section]
\newtheorem{corollary}[theorem]{Corollary}

\newtheorem{lemma}[theorem]{Lemma}
\newtheorem{proposition}[theorem]{Proposition}

\theoremstyle{remark}
\newtheorem{remark}[theorem]{Remark}

\begin{document}
\title{Fixed point data of finite groups\\acting on 3--manifolds} 
\authors{P\'eter E. Frenkel}                  
\coverauthors{P\noexpand\'eter E. Frenkel}                  
\asciiauthors{Peter E. Frenkel}
\address{Department of Geometry, Mathematics Institute\\Budapest 
University of Technology and Economics, 
Egry J. u. 1.\\1111 Budapest, Hungary }                  

\email{frenkelp@renyi.hu}                     

\begin{abstract} 
We consider fully effective ori\-en\-ta\-tion-pre\-serv\-ing smooth
actions of a given finite group $G$ on smooth, closed, oriented 
3--manifolds $M$.
We investigate the relations that necessarily hold between the
numbers of fixed points of various non-cyclic subgroups. In
Section~\ref{arb}, we show that all such relations are in fact equations mod 2, and we
explain how the number of independent equations yields information
concerning low-dimensional equivariant cobordism groups.  Moreover, we restate
 a theorem of A. Sz\H ucs asserting that under the conditions imposed on a
smooth action of $G$ on $M$ as above, the number of $G$-orbits of points $x\in
M$ with non-cyclic stabilizer $G_x$ is even, and we prove the result by using
arguments of G. Moussong.   In Sections~\ref{dih} and
\ref{pla}, 
we determine all the equations for non-cyclic subgroups $G$ of $SO(3)$.
\end{abstract}

\asciiabstract{We consider fully effective orientation-preserving
smooth actions of a given finite group G on smooth, closed, oriented
3-manifolds M.  We investigate the relations that necessarily hold
between the numbers of fixed points of various non-cyclic
subgroups. In Section 2, we show that all such relations are in fact
equations mod 2, and we explain how the number of independent
equations yields information concerning low-dimensional equivariant
cobordism groups.  Moreover, we restate a theorem of A. Szucs
asserting that under the conditions imposed on a smooth action of G
on M as above, the number of G-orbits of points x in M with
non-cyclic stabilizer G_x is even, and we prove the result by using
arguments of G. Moussong.  In Sections 3 and 4, we
determine all the equations for non-cyclic subgroups G of SO(3).}

\asciikeywords{3-manifold, group action, fixed points, equivariant cobordism }

\maketitle

\section{Preliminaries}\label{prel}

We shall use the word ``representation'' to mean a representation
of a group by means of proper rotations of a three-dimensional
Euclidean linear space.

Let $H$ be a subgroup of the group $G$, and let $\rho$ be an
equivalence class of representations of $H$. Denote by
$N_G(H,\rho)$ the subgroup of $G$ formed by those elements $g\in
G$ whose conjugation action fixes $(H,\rho)$ (that is, $g$
normalizes $H$ and the conjugation action leaves the equivalence
class $\rho$ invariant).

We shall use the word ``manifold'' to mean  a smooth, closed,
oriented three-dimensional manifold (if not otherwise stated). All
group actions shall  be assumed smooth and
ori\-en\-ta\-tion-preserving. Consider an action $G\to {\rm
{Diff}}^+( M)$ of the finite group $G$ on the manifold $M$.
Suppose that no nontrivial element of $G$ acts via the identity map on any
component of $M$. (We shall refer to this assumption by saying
that the action of $G$ is ``fully effective''. If $M$ is
connected, then ``fully effective'' is the same as ``effective''.)

Choose a $G$--invariant Riemannian metric on $M$. Then the fixed
point set $M^g$ of an arbitrary element $1\neq g\in G$ is a finite
disjoint union of closed geodesics without self-intersections. For
$x\in M^g$, the proper orthogonal transformation $T_xg\in
SO(T_xM)$ is a rotation of finite order about the axis $T_xM^g$.
It follows that if $M^g$ and $M^h$ have a common component, then
$g$ and $h$ generate a cyclic subgroup of $G$. Therefore, the
fixed point set $M^H$ of any non-cyclic subgroup $H$ of $G$ is
finite. So is the set $Q$ of all points with non-cyclic
stabilizer.

For a non-cyclic subgroup $H$ of $G$, denote by $n^H$ the number
of $N_G(H)$--orbits of
points $x\in M$ whose stabilizer $G_x$ is exactly $H$. 
This is the number of such points divided by the index
$|N_G(H):H|$. In  other words, $n^H$ is the number of $G$--orbits
of points $x\in M$ whose stabilizer $G_x$ is a conjugate of $H$.

 For any subgroup $H$ of $G$ and equivalence class $\rho$ of faithful representations of
$H$, denote by $M^H_\rho$ the set of points $x\in M^H$ such that
the induced representation of $H$ on $T_xM$ is of the equivalence
class $\rho$.
 For non-cyclic $H$, denote by $n^H_\rho$ the number of $N_G(H,\rho)$--orbits of
points $x\in M^H_\rho$ such that $G_x=H$ (these shall be called
``points of type $(H,\rho)$"). This is the number of such points
divided by the index $|N_G(H,\rho):H|$. In  other words,
$n^H_\rho$ is the number of $G$--orbits of points $x\in M$ whose
stabilizer $G_x$ together with its induced representation on
$T_xM$  is a conjugate of $(H,\rho)$.

 We obviously have $n^H=\sum_\rho n^H_\rho$, where $\rho$ runs over those
faithful representations of $H$ that are inequivalent even if
conjugation by elements of $N_G(H)$ is allowed. We also have
$n^H_\rho= n^{H'}_{\rho'}$ if $(H,\rho)$ and ${(H',\rho')}$ are
equivalent under the conjugation action of $G$.

 Note that the sets
$M^H$ and $M^H_\rho$ and the numbers $n^H$ and $n^H_\rho$ do not
depend on the $G$--invariant metric that we have chosen.

\section{Arbitrary groups}\label{arb}

For any finite group $G$, consider the additive Abelian group
$\mathcal A=\mathcal A(G)$ of all integer-valued $G$--invariant
functions
defined on the set of 
 pairs $(H,\rho)$
with $H$ a non-cyclic subgroup of $G$ and $\rho$ an equivalence
class of faithful representations of $H$. We shall be concerned
with the set $\mathcal F_+=\mathcal F_+(G)$ of those functions in
$\mathcal A(G)$ that are realizable as the function $n$ associated
to a fully effective action of $G$ on a manifold. We define
$\mathcal F=\mathcal F(G)$ as the additive subgroup of $\mathcal A
(G)$ generated by $\mathcal F_+ (G)$.

\begin{remark}
Any $n\in\mathcal F_+(G)$ is realizable by an effective action of
$G$ on a connected manifold.
\end{remark}

\begin{proof}
Consider a realization on a manifold $M$. If $M$ is not connected,
then choose points $x$ and $y$ in different components with
trivial stabilizer. Delete small neighborhoods that correspond to
the open 3--ball of radius 1/2 under a
 diffeomorphism of a greater neighborhood onto the open 3--ball of radius 1.
Identify the  boundaries of the deleted neighborhoods in a way
compatible with the given ori\-en\-ta\-tions. Do the same for
$g(x)$ and $g(y)$ $(g\in G)$, using the $g$--image of the
neighborhoods used for $x$ and $y$. The resulting manifold admits
a fully effective $G$--action with the same numbers $n^H_\rho$, and
has  fewer components than $M$ does. Iterating this procedure, we
arrive at a connected manifold.
\end{proof}

\begin{lemma}\label{Absubgr}
{\rm(i)}\qua $\mathcal F_+$ is  the set of non-negative functions in the
additive group $\mathcal F$.

{\rm(ii)}\qua $\mathcal F\geq 2\mathcal A$.
\end{lemma}

\begin{proof}
{\rm(i)}\qua  The action of $G$ on $G\times M$ via left translations on $G$
shows that $0\in \mathcal F_+$.  Taking disjoint union of
manifolds shows that $\mathcal F_+$ is closed under addition.

Suppose that $n, n'\in\mathcal F_+$ and $ n'\geq  n$. We prove
that $n'-n\in\mathcal F_+$. To this end, choose fully effective
actions of $G$ on manifolds $ M$ and $ M'$ that give rise to $ n$
and $ n'$, respectively. A suitable $G$--invariant neighborhood of
$Q\subset M$ (the set of points with non-cyclic stabilizer),
consisting of small open balls, is $G$--equivariantly diffeomorphic
to a $G$--invariant neighborhood of a $G$--invariant subset of $
Q'\subset M'$. The diffeomorphism can be chosen to be
ori\-en\-ta\-tion-reversing, because central reflection in the
origin of $\mathbb{R}^3$ commutes with any proper rotation of
$\mathbb{R}^3$ and reverses ori\-en\-ta\-tion. Delete both
neighborhoods and identify their boundaries to get a closed
manifold  endowed with a fully effective $G$--action that proves
that $n'- n\in\mathcal F_+$.

\medskip

{\rm(ii)}\qua   We first construct a fully effective action of $G$ on a
manifold $M$ 
that has $n^H_\rho>0$  for an arbitrarily chosen non-cyclic
subgroup $H$ of $G$ and an arbitrarily chosen faithful representation $\rho$ of $H$. To
this end, we let the subgroup $H$ act on
$S^3=\mathbb{R}^3\cup\{\infty\}$ by rotations of $\mathbb{R}^3$
given by the representation $\rho$, and we form the twisted
product $M=G\underset H\times S^3$.

Choose a point $x\in M$ of type $(H,\rho)$. Take two copies of
$M$, and delete the two copies of a neighborhood of $Q\backslash
Gx\subset M$, consisting of small open balls. Identify their
boundaries in a $G$--equivariant and ori\-en\-ta\-tion-reversing
way to get a closed manifold endowed with a fully effective
$G$--action that  has $n^H_\rho=2$, and has $n^{H'}_{\rho'}=0$ if
$(H',\rho')$ is not equivalent to $(H,\rho)$ under the conjugation
action of $G$.
\end{proof}

Andr\'as Sz\H ucs called my attention to the fact that the
quotient groups $\mathcal A/\mathcal F$ and $\mathcal F/2\mathcal
A$ can be interpreted in terms of low-dimensional equivariant
cobordism groups. These interpretations are given in the following
two propositions. We do not use them in the sequel. The proofs are left to the reader.

An action of $G$ on a manifold with boundary is called quasi-free  if
all stabilizers are cyclic subgroups of $G$.

\begin{proposition}\label{quasi-free}
The quotient group $\mathcal A(G)/\mathcal F(G)$ is the quasi-free
cobordism group of equivariantly (but not necessarily
quasi-freely) nullcobordant closed oriented two-dimen\-sion\-al
manifolds endowed with fully effective 
$G$--actions.
\end{proposition}

Note that if $\mathcal A(G)\neq 0$  (that is, if $G$ has a subgroup isomorphic
to a non-cyclic subgroup of $SO(3)$), then $\mathcal A(G)/\mathcal F(G)\neq
0$. This shall be seen in  Corollary~\ref{sum}.

Let $\Omega_3(G)$ stand for the oriented equivariant cobordism
group of closed 3--manifolds with free $G$--actions, and let
$\hat\Omega_3 (G)$ stand for the oriented equivariant cobordism
group of those  with arbitrary fully effective $G$--actions.

\begin{proposition}\label{exact}
The sequence $$\Omega_3(G) \to\hat\Omega_3(G)\to\mathcal
F(G)/2\mathcal A(G)\to 0$$ is exact.
\end{proposition}

I do not know if it is true that the map on the left is always
zero, so that $\hat\Omega_3(G)=\mathcal F(G)/2\mathcal A(G)$
(cf. \cite[Corollary 3.4]{Kh},
\cite{Ko1}, \cite{Ko2}, and \cite{S} for results of this kind in the unoriented case).
At any rate, the dimension of the $\mathbb{F}_2$--linear space 
$\hat\Omega_3(G)/2\hat\Omega_3(G)$ is bounded from below by that of $\mathcal
F(G)/2\mathcal A(G)$.

The following theorem of A. Sz\H ucs imposes a  relation
on the values  of any $n\in\mathcal F(G)$, showing in particular
that $\mathcal F(G)\neq\mathcal A(G)$ if $\mathcal A(G)\neq 0$.
The theorem will prove helpful in Section~\ref{pla} when we deal
with the octahedral group.

\begin{theorem}\label{schw}
Let $G$ act fully effectively on $M$. Then the set $Q$ of those
points $x\in M$ that have a non-cyclic stabilizer consists of an
even number of $G$--orbits.
\end{theorem}

The following proof  is due to G. Moussong.
\begin{proof}
Let $P$ stand for the set of points in $M$ whose stabilizer is
non-trivial. Then $P$ is a graph (one-dimensional complex) whose
set of vertices (zero-dimensional skeleton) is $Q$. The set
$P\subset M$ is $G$--invariant, and if $g\in G$ fixes an interior
point of an edge of $P$, then it fixes each point of that edge.
Therefore, $P/G$ is a graph whose set of vertices is $Q/G$.

Choose any point $x\in Q$. The degree in $P/G$ of the vertex
$Gx\in Q/G$ is the number of $G_x$--orbits of directed axes of the
rotations given by the faithful representation $G_x\rightarrowtail
SO(T_xM)$. This number is always 3, as is well known from the
classification  of finite non-cyclic subgroups of $SO(3)$. So the
graph $P/G$ is 3--regular, and therefore it has an even number of
vertices.
\end{proof}

\begin{corollary}\label{sum}
Let $n\in \mathcal F(G)$. Then the sum $\sum n^H_\rho$, with
$(H,\rho)$ running over all $G$--inequivalent pairs, is even.
\end{corollary}

 Note that if $n^G\neq
0$ for some $n\in \mathcal F(G)$, then $G$ is (isomorphic to) a
finite non-cyclic subgroup of $SO(3)$. We shall explicitly
describe $\mathcal F(G)$ (or, equivalently, $\mathcal F_+(G)$) for
these groups $G$ in Section~\ref{dih} and Section~\ref{pla}.
T. tom Dieck \cite[Theorem 6]{D} solved a similar problem for the group
$Z_2^k$ acting on an unoriented manifold of arbitrary dimension. That
investigation 
deals only with the fixed points of the entire group though.

\section{Dihedral groups}\label{dih}
In this section, we explicitly describe $\mathcal F$ for the
dihedral group $D_m$.

For $m>2$, there are exactly $\varphi (m)/2$ inequivalent faithful
representations of the cyclic group $Z_m$. Each extends in a
unique way to a faithful representation of the dihedral group
$D_m$. The two-element group $Z_2$ and the  Kleinian group $D_2$
have unique faithful representations (up to equivalence).

The subgroups of the dihedral group $D_m$ are dihedral groups
themselves. For each $k$ dividing $m$, there are $m/k$ subgroups
of isomorphism type $D_k$. If $m/k $ is odd, then each one is its
own normalizer and they are all conjugates of each other. If $m/k$
is even, then each one has a $D_{2k}$ as its normalizer and they
form two classes that shall be marked by $\bigtriangleup$ and
$\bigtriangledown$.  For any representatives of the two classes,
$D_k^\bigtriangleup\cap D_k^\bigtriangledown=Z_k$.
So a function $n\in \mathcal A(D_m)$ is given by its values
$n^{D_k}_\rho$ (for $m/k$ odd and $\rho$ a faithful representation
of $Z_k$), and $n^{D_k^\bigtriangleup}_\rho$ and
$n^{D_k^\bigtriangledown}_\rho$ (for $m/k$ even and $\rho$ a
faithful representation of $Z_k$).

\begin{theorem}\label{D}
Let $m\geq 2$. The function $n\in\mathcal A(D_m)$ is in $\mathcal
F(D_m)$ if and only if the following two conditions hold for all
$k\geq 2$ dividing $m$.

{\rm(i)}\qua If $m/k$ is odd, then $n^{D_k}_\rho$ is even for all $\rho$.

{\rm(ii)}\qua If $m/k$ is even, then $n^{D_k^\bigtriangleup}_\rho$ and
$n^{D_k^\bigtriangledown}_\rho$ have the same parity for all
$\rho$.
\end{theorem}

\begin{proof} First suppose that $n\in\mathcal F(D_m)$. We prove that (i) and (ii) hold. We may
assume that $n\in\mathcal F_+(D_m)$. Let $n$ be realized by a
fully effective action on $M$.

{\rm(i)}\qua The claim is equivalent to saying that there is an even number
of points of type $(D_k,\rho)$ for any subgroup  of isomorphism
type $D_k$ and for any $\rho$. We shall prove that the number of
such points on any  component $S^1$ of $M^{Z_k}_\rho$ is either
two or zero. Suppose that the point $x\in S^1$ is of type
$(D_k,\rho)$. Consider the action of the two-element quotient
group $D_k/Z_k$
on  $S^1$.
The non-trivial element fixes $x$ and reverses ori\-en\-ta\-tion
on $S^1$. Therefore, it  has a unique fixed point on $S^1$
besides $x$, which is also of type $(D_k,\rho)$. %Then $y\in M^D_m$ and the induced

{\rm(ii)}\qua Suppose to the contrary that, say,
$n^{D_k^\bigtriangleup}_\rho$ is odd and
$n^{D_k^\bigtriangledown}_\rho$ is even. Using
Lemma~\ref{Absubgr}, we may assume that
$n^{D_k^\bigtriangleup}_\rho=1$ and
$n^{D_k^\bigtriangledown}_\rho=0$. This means that, for any
representatives of the two classes, there are no points of type
$(D_k^\bigtriangledown,\rho)$ and there are exactly two points,
say $x$ and $y$, of type $(D_k^\bigtriangleup,\rho)$. The argument
that we have used to prove (i) shows that the points $x$ and $y$
must be on the same component $S^1$ of $M^{Z_k}_\rho\backslash
M^{Z_{2k}}$. For the $D_{2k}$ subgroup generated by the chosen
$D_k^\bigtriangleup$ and $ D_k^\bigtriangledown$ subgroups, the
coset $D_{2k}\backslash D_k^\bigtriangleup$ must interchange $x$
and $y$. So $D_{2k}$ must leave $S^1$ invariant. The Kleinian
quotient group $D_{2k}/Z_k$ acts on $S^1$. The coset
$Z_{2k}\backslash Z_k$ does not have fixed points on $S^1$, so it
preserves ori\-en\-ta\-tion on $S^1$. The subgroup
$D_k^\bigtriangleup$ does have fixed points on $S^1$, so the coset
$D_k^\bigtriangleup\backslash Z_k$ reverses ori\-en\-ta\-tion on
$S^1$. Therefore, the coset $$D_k^\bigtriangledown\backslash
Z_k=(Z_{2k}\backslash Z_k)(D_k^\bigtriangleup\backslash Z_k)$$
also reverses ori\-en\-ta\-tion
 and must have two fixed points. These must be of type $(
D_k^\bigtriangledown,\rho)$, which is a contradiction.

To prove the converse, choose any $k$ dividing $m$, and choose a
faithful representation $\rho$ of $Z_k$.  
Pick a generator $g$ in $Z_m$ for which the generator $g^{m/k}$ of
$Z_k$ is represented by a rotation through $2\pi/k$. Let $g$ act
on $S^3\subset\mathbb{C}^2$ by $$g(z,w)=(e^{2\pi ik/m}z, e^{2\pi
i/m}w).$$ Choose an arbitrary element of the coset $D_m\backslash
Z_m$ and let it act by complex conjugation on both coordinates.
This defines an action of $D_m$ on $S^3$. The set of points whose
stabilizer non-trivially intersects $Z_m$ is the circle
$$S^1=\{w=0\}=M^{Z_k}_\rho.$$ All points $x$ of $S^1$ have
$G_x\cap Z_m=Z_k$, so $G_x$ is either $Z_k$ or a $D_k$ subgroup.
The two-element quotient group $D_k/Z_k$ of each $D_k$ subgroup
acts on $S^1$. The coset $D_k\backslash Z_k$  reverses
ori\-en\-ta\-tion on $S^1$ and therefore has exactly two fixed
points. If $m/k$ is even, then this means that
$n^{D_k^\bigtriangleup}_\rho= n^{D_k^\bigtriangledown}_\rho=1$ and
all other numbers $n^H_\rho$ are zero.

The theorem now follows from Lemma~\ref{Absubgr}.
\end{proof}

\begin{corollary}
{\rm(i)}\qua The dimension of the $\mathbb{F}_2$--linear space $\mathcal
A(D_m)/\mathcal F(D_m)$ is $\lfloor m/2\rfloor$.

{\rm(ii)}\qua The dimension of the $\mathbb{F}_2$--linear space $\mathcal
F(D_m)/2\mathcal A(D_m)$ is $\lfloor m/4\rfloor$ if $m$ is even
and is zero if $m$ is odd.
\end{corollary}

\begin{proof}
{\rm(i)}\qua  The dimension is the total number of inequivalent faithful
representations of the cyclic groups $Z_k$ with $k\geq 2$ dividing
$m$, which is $$\sum_{2<k\mid m} \varphi (k)/2 \qquad{\rm { resp.
}}\qquad 1+\sum_{2<k\mid m} \varphi (k)/2 $$ as $m$ is odd or
even.  The claim follows from the fact that $\sum_{k\mid m}\varphi
(k)/2=m/2$ and $\varphi (1)=\varphi (2)=1$.

{\rm(ii)}\qua The dimension is the total number of inequivalent faithful
representations of the cyclic groups $Z_k$ with $k\geq 2$ and
$m/k$ even. This equals zero for odd $m$ and equals
$\lfloor(m/2)/2\rfloor=\lfloor m/4\rfloor$ for even $m$.
\end{proof}

\section{Groups of rotations of the Platonic solids}\label{pla}
In this section, we explicitly describe  $\mathcal F$ for the
tetrahedral, the octahedral, and the icosahedral group. We shall
frequently use the expression ``a conjugation action of $G$ on
$SO(3)$", which shall mean
 an injective homomorphism $G\rightarrowtail SO(3)$ composed with the
conjugation action $SO(3)\rightarrowtail {\rm Diff}^+(SO(3))$ of
$SO(3)$ on $SO(3)$. Note that a  suitable diffeomorphism between
$SO(3)$ and $\mathbb{R}P^3$ turns this into an action of $G$ on
$\mathbb{R}P^3=\mathbb{R}^3\cup\mathbb{R}P^2$ by rotations of
$\mathbb{R}^3$ extended to $\mathbb{R}P^3$.

The only non-cyclic proper subgroup of the tetrahedral group $A_4$
is the Kleinian group $D_2$. Faithful representations of both
$A_4$ and $D_2$ are unique up to equivalence. So a function
$n\in\mathcal A(A_4)$ is given by its two values $n^{A_4}$ and
$n^{D_2}$.

\begin{theorem}\label{A4}
The function $n\in\mathcal A(A_4)$ is in $\mathcal F (A_4)$ if and
only if $n^{D_2}$ and $n^{A_4}$ are of the same parity.
\end{theorem}

\begin{proof}
 ``Only if" is a particular case of Corollary~\ref{sum}. It also follows from 
Theorem~\ref{D}(i) when applied to $k=m=2$.

To prove the converse, observe that $n^{A_4}
=1$ holds for the conjugation action of $A_4$ on $SO(3)$. The
theorem now follows from Lemma~\ref{Absubgr}.
\end{proof}

\begin{corollary}
{\rm(i)}\qua $\dim_{\mathbb{F}_2}\mathcal A(A_4)/\mathcal
 F(A_4)=1.$

{\rm(ii)}\qua $\dim_{\mathbb{F}_2}\mathcal F(A_4)/2\mathcal A(A_4)=1.$
\end{corollary}

 The non-cyclic  subgroups of the octahedral group $S_4$ are of isomorphism type
$S_4$, $A_4$, $D_4$, $D_3$ and $D_2$. All have unique faithful
representations. Any two of them that are isomorphic are
conjugates of each other, except for one of the four Kleinians
which is normal in $S_4$ and shall be denoted by $D_2^*$. So a
function $n\in\mathcal A(S_4)$ is given by its six values
$n^{S_4}$, $n^{A_4}$, $n^{D_4}$, $n^{D_3}$, $n^{D_2}$ and
$n^{D_2^*}$.

\begin{theorem}\label{S4}
The function $n\in\mathcal A(S_4)$ is in $\mathcal F(A_4)$ if and
only if the four numbers $n^{D_3}$, $n^{D_4}$, $n^{S_4}$ and
$n^{D_2^*}+n^{D_2}+n^{A_4}$ are of the same parity.
\end{theorem}

\begin{proof}
To prove ``only if", suppose that  $n\in\mathcal F_+(S_4)$ is
realized
  by an action of $S_4$ on a manifold. Restriction of that action to
 a  $D_3$  subgroup
 shows that
$n^{D_3}+n^{S_4}\in\mathcal F_+(D_3)$. Restriction to a $D_4$
subgroup shows that $n^{D_4}+n^{S_4}$ is the value at $D_4$ of a
function in $\mathcal F (D_4)$.  Theorem~\ref{D}(i), when applied
to $k=m=3$ and to $k=m=4$, shows that $n^{D_3}$, $n^{D_4}$ and
$n^{S_4}$ are of the same parity. Corollary~\ref{sum} says that
$$n^{D_3}+n^{D_4}+n^{S_4}+n^{D_2^*}+n^{D_2}+n^{A_4}$$ is even, so
``only if" is proved.

To prove the converse, consider the conjugation actions of the
subgroups $H=S_4$, $A_4$ and $D_4$ on $SO(3)$, and form twisted
products $S_4\underset H\times SO(3)$. The functions $n$
associated to these three actions of $S_4$ are linearly
independent mod 2, since they take the  values listed below at
$S_4$, $A_4$ and $D_2$. $$\begin{matrix}  &  SO(3)  &
{S_4\underset{A_4}\times SO(3)}  &  {S_4\underset{D_4}\times
SO(3)}
\\
{S_4}  &  1  &  0  &  0\\ {A_4}  &  *  &  1  &  0\\ {D_2}  &  *  &
*  &  1\end{matrix}$$ The theorem now follows from
Lemma~\ref{Absubgr}.
\end{proof}

\begin{corollary}
{\rm(i)}\qua $\dim_{\mathbb{F}_2}\mathcal A(S_4)/\mathcal F(S_4)=3$.

{\rm(ii)}\qua $\dim_{\mathbb{F}_2}\mathcal F(S_4)/2\mathcal A(S_4)=3$.
\end{corollary}

 The non-cyclic proper subgroups of the icosahedral group $A_5$ are of 
isomorphism type
$A_4$, $D_5$, $D_3$ and $D_2$. Any two of them that are isomorphic
are conjugates of each other. The icosahedral group has exactly
two inequivalent faithful representations. We shall denote them
by $+$ and $-$. Their restriction to any of the six $D_5$
subgroups yields the two ($=\varphi (5)/2$) inequivalent faithful
representations of that subgroup. These shall also be denoted by
$+$ and $-$, respectively. The faithful representations of the
other non-cyclic subgroups are unique up to equivalence. So a
function $n\in\mathcal A(A_5)$ is given by its seven values
$n^{A_5}_+$, $n^{A_5}_-$, $n^{A_4}$, $n^{D_5}_+$, $n^{D_5}_-$,
$n^{D_3}$ and $n^{D_2}$.

\begin{theorem}\label{A5}
The function $n\in\mathcal A( A_5)$ is in $\mathcal F(A_5)$ if and
only if the four sums $$\begin{aligned} n^{D_2}+n^{A_4}  &  &  & &
+n^{A_5}_+ & +n^{A_5}_- & ,\\  & n^{D_3} &
 &  & +n^{A_5}_+ & +n^{A_5}_- & ,\\  &  & n^{D_5}_+ &  & +n^{A_5}_+ &  & ,\\ &  &  &
n^{D_5}_- &  & +n^{A_5}_- & \end{aligned}$$ are even.
\end{theorem}

\begin{proof}
``Only if" follows from Theorem~\ref{D}(i) when applied to
$k=m=2$, $3$, $5$ and $5$, respectively.

To prove the converse, consider the two conjugation actions of
$A_5$ on $SO(3)$. Also consider the conjugation action of an $A_4$
subgroup on $SO(3)$, and form the twisted product $A_5\underset
{A_4}\times SO(3)$.  The functions $n$ associated to these three
actions of $A_5$ are linearly independent $\mod 2$, since they
take the values listed below at $(A_5,+)$, $(A_5,-)$ and $A_4$.
$$\begin{matrix}   &  { SO(3)}  &  {SO(3)}   &
{A_5\underset{A_4}\times SO(3)}
\\
{(A_5,+)}   &  1  &  0 &  0\\ {(A_5,-)}   &  0  &  1 &  0\\ {A_4}
 &  0  &  0 &  1\end{matrix}$$ The theorem now follows from
Lemma~\ref{Absubgr}.
\end{proof}

\begin{corollary}
{\rm(i)}\qua $\dim_{\mathbb{F}_2}\mathcal A(A_5)/\mathcal F(A_5)=4$.

{\rm(ii)}\qua $\dim_{\mathbb{F}_2}\mathcal F(A_5)/2\mathcal A(A_5)=3$.
\end{corollary}

\subsection*{Acknowledgements}
I am deeply indebted to Professor Andr\'as Sz\H ucs for calling my
attention to the problem, and for helpful discussions.
I am also grateful to Professor Andr\'as Szenes for useful comments.

This  research was partially supported by OTKA grant  T--042769.

\Addressesr

\end{document}

%% file: agtout.tex
\def\ifplaintex{\expandafter\ifx\csname documentclass\endcsname\relax}

\def\gtp{{\mathsurround=0pt\it $\cal G\mskip-2mu$eometry \&\ 
$\cal T\!\!$opology $\cal P\!$ublications}}  % GT publications

\def\Addressesr{\bigskip
{\small \parskip 0pt \leftskip 0pt \rightskip 0pt plus 1fil \def\\{\par}
\sl\theaddress\par
\medskip
\rm Email:\stdspace\tt\theemail\hfill\rm Received:\qua\receiveddate \par}}

\def\recd{{\small Received:\qua\receiveddate\ifx\reviseddate\relax
\else\qquad Revised:\qua\reviseddate\fi\par}} 

\def\lognumber#1{\def\thelognumber{#1}}
\def\volumenumber#1{\def\thevolumenumber{#1}}
\def\volumeyear#1{\def\thevolumeyear{#1}}
\def\papernumber#1{\def\thepapernumber{#1}}
\def\pagenumbers#1#2{\def\startpage{#1}\def\finishpage{#2}}
\def\published#1{\def\publishdate{#1}}

\def\received#1{\def\receiveddate{#1}}

\def\accepted#1{\def\accepteddate{#1}}

\def\asciiauthors#1{\def\theasciiauthors{#1}}

\def\coverauthors#1{\def\thecoverauthors{#1}}
\long\def\asciiabstract#1{\long\def\theasciiabstract{#1}}
\def\asciikeywords#1{\def\theasciikeywords{#1}}

\let\\\par\let\thelognumber\relax\let\thevolumenumber\relax
\let\thepapernumber\relax\let\thevolumeyear\relax\let\startpage\relax
\let\finishpage\relax\let\publishdate\relax\let\receiveddate\relax
\let\reviseddate\relax\let\accepteddate\relax\let\theasciititle\relax
\let\theasciiauthors\relax
\let\theasciiabstract\relax\let\theasciikeywords\relax

\let\thecoverauthors\relax\let\theasciiemail\relax

\ifplaintex
\font\logobig=cmssbx10 scaled 3836
\font\logomed=cmssbx10 scaled 2557
\else
\font\logobig=cmssbx10 scaled 4200
\font\logomed=cmssbx10 scaled 2800
\fi

\long\def\makeagttitle{   %%% start of definition of \makeagttitle
\count0=\startpage
\agt\hfill      %   Journal title (top left) 
\hbox to 45truept{\vbox to 0pt{\vglue -13truept{\logomed A\kern -.37em{\logobig 
T}\kern -.38em G}\vss}\hss}
\break
{\small Volume \thevolumenumber\ (\thevolumeyear)
\startpage--\finishpage\nl
Published: \publishdate}

\vglue .25truein

{\parskip=0pt\leftskip 0pt plus
1fil\def\\{\par\smallskip}{\Large\bf\thetitle}\par\medskip} \vglue
0.05truein

{\parskip=0pt\leftskip 0pt plus 1fil\def\\{\par}{\sc\theauthors}
\par\medskip}%
 
\vglue 0.03truein 

{\small\leftskip 25truept\rightskip 25truept{\bf Abstract}\stdspace\theabstract

{\bf AMS Classification}\stdspace\theprimaryclass
\ifx\thesecondaryclass\relax\else; \thesecondaryclass\fi\par
{\bf Keywords}\stdspace \thekeywords\par}\vglue 7truept

}   %%%% end of definition of \makeagttitle

\ifplaintex
\font\phead=cmsl9 scaled 950
\font\pnum=cmbx10 scaled 913
\font\pfoot=cmsl9 scaled 950
\headline{\vbox to 0pt{\vskip -4.5mm\line{\small\phead\ifnum
\count0=\startpage ISSN 1472-2739 (on-line) 1472-2747 (printed)
\hfill {\pnum\folio}\else\ifodd\count0\def\\{ }% 
\ifx\theshorttitle\relax\thetitle\else\theshorttitle\fi\hfill{\pnum\folio}
\else\def\\{ and }{\pnum\folio}\hfill\ifx\theshortauthors\relax\theauthors
\else\theshortauthors\fi\fi\fi}\vss}}
\footline{\vbox to 0pt{\vglue 0mm\line{\small\pfoot\ifnum\count0=\startpage
\copyright\ \gtp\hfill\else
\agt, Volume \thevolumenumber\ (\thevolumeyear)\hfill\fi}\vss}}
\else
\font=cmsl9 scaled 1050
\font\lnum=cmbx10 
\font=cmsl9 scaled 1050
\makeatletter
\def\@oddhead{{\small\lhead\ifnum\count0=\startpage ISSN 1472-2739 
(on-line) 1472-2747 (printed)\hfill {\lnum\number\count0}\else\ifodd\count0
\def\\{ }\ifx\theshorttitle\relax \thetitle \else\theshorttitle\fi\hfill
{\lnum\number\count0}\else\def\\{ and }{\lnum\number\count0}
\hfill\ifx\theshortauthors\relax 
\theauthors\else\theshortauthors\fi\fi\fi}}\def\@evenhead{\@oddhead}
\def\@oddfoot{\small\lfoot\ifnum\count0=\startpage\copyright\ \gtp\hfill\else
\agt, Volume \thevolumenumber\ (\thevolumeyear)\hfill\fi}
\def\@evenfoot{\@oddfoot}
\makeatother
\fi
\let\maketitlepage\makeagttitle
\let\makeshorttitle\maketitlepage
\let\maketitle\maketitlepage

\newwrite\gtoutfile
\long\gdef\makeheadfile{  %%% start of definition of \makeheadfile
{\def\\{, }\def\s{ }
\immediate\openout\gtoutfile head.xxx
\immediate\write\gtoutfile{To: math@arxiv.org}
\immediate\write\gtoutfile{Subject: put OR rep NNNNN:ppppp}
\immediate\write\gtoutfile{--text follows this line--}
\immediate\write\gtoutfile{Proxy-for: \ifx\theasciiauthors\relax
\theauthors\else\theasciiauthors\fi\s<\ifx\theasciiemail\relax\theemail\else\theasciiemail\fi>}
\immediate\write\gtoutfile{\noexpand\\}
\immediate\write\gtoutfile{Authors: \ifx\theasciiauthors\relax
\theauthors\else\theasciiauthors\fi}
{\def\\{ }\immediate\write\gtoutfile{Title: \ifx\theasciititle\relax
\thetitle\else\theasciititle\fi}}
\immediate\write\gtoutfile{Subj-class: GT or SG, GR etc}
\immediate\write\gtoutfile{MSC-class: \theprimaryclass\ifx\thesecondaryclass\relax\else, \thesecondaryclass\fi}
\immediate\write\gtoutfile{Journal-ref: Algebr. Geom. Topol. \thevolumenumber\s
(\thevolumeyear) \startpage-\finishpage}
\immediate\write\gtoutfile{Comments: Published by Algebraic and
Geometric Topology at}
\immediate\write\gtoutfile{\s\s\s  http://www.maths.warwick.ac.uk/agt/AGTVol\thevolumenumber/agt-\thevolumenumber-\thepapernumber.abs.html}
\immediate\write\gtoutfile{\noexpand\\}
\immediate\write\gtoutfile{}
\ifx\theasciiabstract\relax
\immediate\write\gtoutfile{\theabstract}\else
\immediate\write\gtoutfile{\theasciiabstract}\fi
\immediate\write\gtoutfile{}
\immediate\write\gtoutfile{\noexpand\\}
\immediate\write\gtoutfile{}
\immediate\closeout\gtoutfile}}  %%% end of definition of \makeheadfile

\def\maketitlepage{\makeagttitle\makeheadfile}
\let\makeshorttitle\maketitlepage
\let\maketitle\maketitlepage